\documentclass[11pt]{amsart}
\usepackage{graphicx}
\usepackage{latexsym}
\usepackage{amsfonts,amsmath,amssymb}
\newtheorem{theorem}{Theorem}[section]																
\newtheorem{lemma}[theorem]{Lemma}

\newtheorem{corollary}[theorem]{Corollary}

\usepackage{color}

\def\a{\alpha}

\def\ga{\gamma}
\def\part{\partial}

\def\b1{\bold 1}

\newcommand{\beq}{\begin{equation}}
\newcommand{\eeq}{\end{equation}}

\theoremstyle{remark}

\numberwithin{equation}{section}
\date{\today}

\begin{document}

\title[Branching process]{Counting subtrees of the branching process tree by the number of leaves} 
\date{}
\author{Boris Pittel}
\address{Department of Mathematics, Ohio State University, 231 West 18-th Avenue, Columbus OH 43210-1175}
\email{pittel.1@osu.edu}
\begin{abstract} We study the distribution of the number of leaves of the subtree chosen uniformly at random among all the subtrees of the critical branching process tree at extinction.

\end{abstract}
\keywords
{Branching process, random tree, extinction, asymptotics}

\subjclass{60C05; 05C05, 92B10}

\maketitle
\section{Introduction and results} 
Consider a branching process initiated by a single progenitor. This process is visualized as a growing rooted tree.
The root is the progenitor, connected by edges to each of its immediate descendants (children), that are {\it ordered\/},
say by seniority.
Each of the children becomes the root of the corresponding subtree, so that the children of all these roots are the grandchildren of the progenitor. We obviously get a recursively defined process. It delivers 
a nested sequence of trees, which is 
either infinite, or terminates at a moment when none of the current leaves have children.

The classic Galton-Watson branching process is the case when the number of each member's children {\bf (a)\/} is independent
of those numbers for all members from the preceding and current generations and {\bf (b)\/} has the same distribution $\{p_j\}_{j\ge 0}$, $(\sum_j p_j=1)$. It is well-known that if $p_0>0$ and $\sum_{j\ge 0}jp_j=1$, then the process terminates with probability $1$, Harris \cite{Har}. A standard example is $p_0=p_2=\tfrac{1}{2}$, in which case we have a binary tree.

Let $T$ denote the terminal tree. Given a finite rooted tree $\mathcal T$, we have
\[
\Bbb P(T=\mathcal T)=\prod_{v\in \mathcal V(\mathcal T)} p_{d(v,\mathcal T)},
\]
where $\mathcal V(\mathcal T)$ is the vertex set of $\mathcal T$, and $d(v, \mathcal T)$ is the out-degree of vertex $v\in V(\mathcal T)$.  Introduce $L=L_T$ the total number of leaves of $T$,
i.e. vertices of $T$ with out-degree $0$, and $X(t)=X_T(t)$ the total number of the (full) subtrees of $T$ with $t$ leaves. So,
the total number of subtrees $\sum_{t\ge 1}X(t)=V(T)=|\mathcal V(T)|$. Introduce $\mathcal L=\mathcal L(T)$, the number of leaves in the subtree chosen uniformly at random among all $V=V(T)$ subtrees of $T$; so, conditionally on $T$, we have $\Bbb P(\mathcal L=t|T)=\tfrac{X(t)}{V}$, and--conditionally on $L$ only--$\Bbb P(\mathcal L=t|L)=\Bbb E\bigl[\tfrac{X(t)}{V}|L\bigr]$. 
It was proved in \cite{Pit1} that $f(x):=\Bbb E[x^L]=\sum_{k\ge 1}x^k \Bbb P(L=k)$ satisfies
\begin{equation}\label{1}
f(x)=\sum_{j\ge 1}p_jf^j(x) +p_0x,\quad |x|\le 1.
\end{equation}
Also $\Bbb P(L<\infty)=1$, and assuming that $\{p_j\}$ has a finite variance $\sigma^2$, 
\begin{equation}\label{1.05}
\Bbb P(L=\ell)=\tfrac{\ga (2\ell-1)!!}{2^{\ell}\ell!}+O(\ell^{-2})=\tfrac{\ga}{\ell^{3/2}}+O(\ell^{-2}),\quad \ga:=\bigl(\tfrac{p_0}{2\pi\sigma^2}\bigr)^{1/2},
\end{equation}
$\sigma^2$ being the variance of $\{p_j\}$. (To be sure, we were interested in the case $p_1=0$, but \eqref{1.05}
holds for all $p_1<1$ as well.) A less sharp  formula $\Bbb P(L=\ell)=\tfrac{\ga}{\ell^{3/2}}(1+o(1))$ has long been
known, Kolchin \cite{Kol3} (Ch. 2, Lemma 4).
\begin{theorem}\label{thm1}
{\bf (i)\/} For the binary case $p_0=p_2=1/2$, 
\begin{equation*}
\Bbb P(\mathcal L=t|L=\ell)=\tfrac{t}{\ell(2\ell-1)}\cdot\tfrac{\binom{\ell}{t}^2}{\binom{2(\ell-1)}{2(t-1)}}.
\end{equation*}
{\bf (ii)\/} Consequently, for $t=o(\ell)$ we have
\begin{equation*}
\Bbb P(\mathcal L=t|L=\ell)=\tfrac{1+O(\ell^{-1/2}+t/\ell)}{2^{2t-1}t}\binom{2(t-1)}{t-1}\!=\!(1+O(\ell^{-1/2}+t/\ell))\Bbb P(L=t),
\end{equation*}
implying tightness of the sequence of distributions $\{\Bbb P(\mathcal L=t|L=\ell)\}_{t\ge 1}$. 
However, $\Bbb E[\mathcal L|L=\ell]\sim\tfrac{\sqrt{\pi\ell}}{2}$, as $\ell\to\infty$.
\end{theorem}
\noindent {\bf Note.\/} So, for large $\ell$, the number of leaves in the  subtree chosen at random from all $(2\ell-1)$ subtrees of the extinction tree with $\ell$ leaves, is distributed asymptotically as the number of leaves in the extinction tree.

\noindent Our second result is for a general critical distribution $\{p_j\}_{j\ge 0}$.
\begin{theorem}\label{thm2} {\bf (i)\/} For each fixed $t\ge 1$,
\[
\lim_{\ell\to\infty}\Bbb P(\mathcal L=t|L=\ell)=(1-p_1)\Bbb P(L=t).
\]
Since $\sum_{t\ge 1}\Bbb P(L=t)=1$, we see that the sequence of the distributions $\{\Bbb P(\mathcal L=t|L=\ell)\}_{1\le t\le \ell}$, $(\ell\ge 1)$, is not
tight, unless (like in the binary case) $\,p_1=0$. In words, $p_1$ is the limiting deficit of the leaf-set size distribution for the random subtree. {\bf (ii)\/} More precisely,
\[
\Bbb P\Bigl(\mathcal L>\tfrac{\ell^{1/2}}{\log^2\ell}\Big| L=\ell\Bigr)=p_1+O(\ell^{-1/4}\log\ell),
\]
so that, conditioned on $\{L=\ell\}$,  $\mathcal L$ exceeds $\tfrac{\ell^{1/2}}{\log^2\ell}$ with conditional probability bounded away from zero as $\ell\to\infty$.
\end{theorem}
\noindent {\bf Note.\/}  Since $p_1<1$, the equation \eqref{1} is equivalent to
\[
f(x)=\sum_{j\ge 2}p_j' f^j(x)+p_0'x,\quad p_j'=\tfrac{p_j}{1-p_1},\,\,j=0,2,3,\dots,
\]
and $p_0', p_2',p_3',\dots$ is a probability distribution, with $\sum_j j p_j'=1$ again, but with $p_1'=0$.  Therefore $L$ and $L'$, the number of leaves in the terminal tree $T'$ associated with $\{p_j'\}$, are {\it equidistributed\/}. So, the part {\bf(i)\/}
can be interpreted as saying that, conditioned on the event ``no father has a single child'',  the uniformly random subtree of $T$, in the limit, has the number of leaves distributed as that for the tree $T'$.

As a source for our inspiration, we should mention the Russian mathematician Valentin Kolchin who pioneered and championed  the study of connection between conditiional branching processes and combinatorics of random trees since the mid-seventies, Kolchin \cite{Kol1}, \cite{Kol2}, and \cite{Kol3}. We refer the reader to David Aldous \cite{Ald} and
Svante Janson \cite{Jan} for two, $20$ years apart, encyclopedic surveys of limit results for the conditioned and the simply generated trees, without convergence rates, that analyze a rich variety of fringe distributions. The distribution of random subtree size is listed in \cite{Ald}  as one of the basic problems in this class of distributions. 

\section{Generating functions identities} Here are two useful identities that follow from implicit differentiation of 
\eqref{1} for $f(x)$:
\begin{equation}\label{1.1}
\sum_{j\ge 1}jp_j f^{j-1}(x)=1-\tfrac{p_0}{f'(x)},\quad \sum_{j\ge 1}j(j-1)p_jf^{j-2}(x)=\tfrac{p_0f^{''}(x)}{(f'(x)^3}.
\end{equation}
In particular, since $f(0)=0$, we have $f'(0)=\tfrac{p_0}{1-p_1}$. 

Our first task is to derive an equation for $f(\ell,t)=\Bbb E[X(t)\Bbb I(L=\ell)]$, where $\Bbb I(B)$ stands for the indicator of an event $B$. Of course, $f(\ell,t)=0$ for $\ell<t$, and $f(t,t)=\Bbb P(L=t)$. Consider $\ell>t$.
With probability $p_j$ the root has $j$ children; let $X_i(t)$ denote the number of subtrees with $t$ leaves in the tree rooted at the $i$-th child. Then
\begin{multline*}
f(\ell,t)=\Bbb E\bigl[X(t)\Bbb I(L=\ell)\bigr]=\sum_{j\ge 1} p_j\Bbb E\biggl[\Bigl(\sum_{i\in [j]}X_i(t)\Bigr)\cdot \Bbb I\Bigl(\sum_{i'\in [j]}L_{i'}=\ell\Bigr)\biggl]\\
=\sum_{j\ge 1} jp_j\Bbb E\biggl[X_1(t)\cdot \Bbb I\Bigl(\sum_{i\in [j]}L_{i}=\ell\Bigr)\biggl]\\
= \sum_{j\ge 1} jp_j\sum_{\mu\le \ell}\Bbb E\bigl[X(t)\Bbb I(L=\mu)\bigr]\cdot\Bbb P\Bigl(\sum_{i\in [j]\setminus \{1\}}L_i=\ell-\mu\Bigr)
\\
=\sum_{j\ge 1} jp_j\sum_{\mu\le \ell} f(\mu,t)\cdot[x^{\ell-\mu}] f^{j-1}(x)\\
=\sum_{\mu\le \ell} f(\mu,t)\bigl[x^{\ell-\mu}\bigr]\sum_{j\ge 1}jp_j f^{j-1}(x).
\end{multline*}
So, in combination with \eqref{1.1}, we obtain
\begin{equation}\label{3}
f(\ell,t)=\sum_{\mu\le \ell} f(\mu,t)\cdot\bigl[x^{\ell-\mu}\bigr]\bigl(1-\tfrac{p_0}{f'(x)}\bigr),\quad \ell>t.
\end{equation}
And we remind that $f(t,t)=\Bbb P(L=t)$. The convolution on the RHS of \eqref{3} positively dictates usage of generating functions. For $|y|<1$, using $f'(0)=\tfrac{p_0}{1-p_1}$, we have
\begin{multline*}
\sum_{\ell\ge t}y^{\ell}f(\ell,t)=y^t\,\Bbb P(L=t) - y^t\, \Bbb P(L=t)\cdot [x^0]\bigl(1-\tfrac{p_0}{f'(x)}\bigr)\\
+\sum_{\ell\ge t}y^{\ell}\sum_{\mu=t}^{\ell}f(\mu,t)\cdot[x^{\ell-\mu}]\bigl(1-\tfrac{p_0}{f'(x)}\bigr)\\
=(1-p_1)y^t\, \Bbb P(L=t)+\sum_{\mu\ge t}f(\mu,t)y^{\mu}\sum_{\ell\ge\mu}y^{\ell-\mu}\cdot [x^{\ell-\mu}]
\bigl(1-\tfrac{p_0}{f'(x)}\bigr)\\
=(1-p_1)y^t\,\Bbb P(L=t)+\bigl(1-\tfrac{p_0}{f'(y)}\bigr)\sum_{\mu\ge t}f(\mu,t)y^{\mu},\\
\end{multline*}
implying that
\begin{equation}\label{4}
\sum_{\ell\ge t}y^{\ell}f(\ell,t)=\tfrac{1-p_1}{p_0}y^t f'(y)\,\Bbb P(L=t).
\end{equation}
Recalling that $f(\ell,t)=f(\ell,t)=\Bbb E[X(t)\Bbb I(L=\ell)]$, and using $[y^s] f'(y)=(s+1)\cdot[y^{s+1}]f(y)$, we arrive at
\begin{lemma}\label{lem1}
\begin{align*}
\Bbb E[X(t)|L=\ell]&=\tfrac{1-p_1}{p_0}\tfrac{[y^{\ell-t}]f'(y)\times \Bbb P(L=t)}{\Bbb P(L=\ell)}=\tfrac{1-p_1}{p_0}\,\tfrac{(\ell-t+1)
\Bbb P(L=\ell-t+1)\,\Bbb P(L=t)}{\Bbb P(L=\ell)}.\\
\end{align*}
\end{lemma}
 \section{Proof of Theorem \ref{thm1}}

Notice that $V=2\ell-1$ on the event $\{L=\ell\}$ for the binary case. 
\begin{proof} By Lemma \ref{lem1} (first identity), with a bit of elementary work, we have 
\begin{multline*}
\Bbb P(\mathcal L=t|L=\ell)=\tfrac{1}{2\ell-1}\Bbb E[X(t)|L=\ell]=\tfrac{2}{2\ell-1}\tfrac{[y^{\ell-t}]f'(y)\times [x^t]f(x)}{[x^{\ell}]f(x)}\\
=\tfrac{1}{2\ell-1}\cdot \tfrac{[y^{\ell-t}](1-y)^{-1/2}\times [x^t](-(1-x)^{1/2})}{[x^{\ell}](-(1-x)^{1/2})}
=\tfrac{t}{\ell(2\ell-1)}\cdot\tfrac{\binom{\ell}{t}^2}{\binom{2(\ell-1)}{2(t-1)}}.
\end{multline*}
The asymptotic formula $\Bbb E[\mathcal L|L=\ell]\sim\tfrac{\sqrt{\pi\ell}}{2}$ follows easily from
\[
t\cdot\Bbb P(\mathcal L=t|L=\ell)=\tfrac{\ell}{2\ell-1}\tfrac{\binom{2(t-1)}{t-1}\binom{2(\ell-t)}{\ell-t}}{\binom{2(\ell-1)}{\ell-1}},
\]
the formula $\binom{2m}{m}\sim \tfrac{2^{2m}}{\sqrt{\pi m}}$, $(m\to\infty)$, and $\int_0^1\tfrac{dx}{\sqrt{x(1-x)}}=\pi$. Finally, by \eqref{1.05} and the second identity in Lemma \ref{lem1}, we have
\begin{multline*}
\Bbb P(\mathcal L=t|L=\ell)=\tfrac{1}{2\ell-1}\Bbb E[X(t)|L=\ell]\\
=\tfrac{2(\ell-t+1)(\ell-t+1)^{-3/2}}{(2\ell-1)\ell^{-3/2}}\cdot \bigl (1+O((\ell-t+1)^{-1/2})\bigr)
\Bbb P(L=t)\\
=(1+O(\ell^{-1/2}+t/\ell))\Bbb P(L=t),
\end{multline*}
provided that $t=o(\ell)$.
\end{proof}

\section{Proof of Theorem \ref{thm2}} For a general $\{p_j\}_{j\ge 0}$ with $p_0>0$, $\sum_j jp_j=1$, $V=\sum_sX(s)$ is random on the event $\{L=\ell\}$, i.e. the number of the subtrees of the extinction tree to choose from is random. So, here
$
\Bbb P(\mathcal L=t|L=\ell)=\Bbb E\bigl[\tfrac{X(t)}{V}\big| L=\ell\bigr].
$
The ratios of dependent random variables can be problematic for evaluation, precise or even asymptotic, of their
moments. Fortunately, here we hope that, for $\ell$ large, conditionally on $\{L=\ell\}$, 
 the distribution of $V$ is sharply concentrated around $\Bbb E[V|L=\ell]$. If so, we can expect that for large $\ell$, 
$\Bbb P(\mathcal L=t|L=\ell)\sim \tfrac{\Bbb E[X(t)|L=\ell]}{\Bbb E[V|L=\ell]}$ for a wide range of $t<\ell$.  And the last fraction is perfectly amenable to asymptotic analysis.

To prove concentration, let us introduce $g(y,\ell)=\Bbb E[y^{V}\Bbb I(L=\ell)]$, $(\ell\ge 1)$, so that $g(y,1)=p_0y$, and define 
$G(y,x):=\sum_{\ell\ge 1}x^{\ell}g(y,\ell)$. For $\ell>1$, we have 
\begin{equation*}
g(y,\ell)=y\sum_{j\ge 1}p_j\sum_{\ell_1+\cdots+\ell_j=\ell}\,\prod_{i\in [j]}g(y,\ell_i)
=[x^{\ell}]\,y\sum_{j\ge 1}p_j G^j(y,x).
\end{equation*}
Consequently
\begin{equation}\label{5}
G(y,x)=p_0 yx +y\sum_{j\ge 1}p_j G^j(y,x).
\end{equation}
We use \eqref{5} to get $G'_y(1,x)$ and $G^{''}_y(1,x)$, since  
\begin{equation}\label{6}
\begin{aligned}
G'_y(1,x)&=\Bbb E[V x^L] \Longrightarrow \Bbb E[V\Bbb I(L=\ell)]=[x^{\ell}]G'_y(1,x)
,\\
G^{''}_y(1,x)&=\Bbb E[V(V-1) x^L]\Longrightarrow \Bbb E[V(V-1)\Bbb I(L=\ell)]=[x^{\ell}]G^{''}_y(1,x).
\end{aligned}
\end{equation}
We have 
\begin{align*}
G(1,x)&=\sum_{\ell\ge 1}x^{\ell}g(1,\ell)=\sum_{\ell\ge 1}x^{\ell}\Bbb P(L=\ell)=f(x),\\
G'_y(1,x)&=\biggl(p_0x+\sum_{j\ge 1}p_jG^j(y,x)+y\sum_{j\ge 1}jp_jG^{j-1}(y,x)G'_y(y,x)\Biggr)\Big|_{y=1}\\
&=f(x)+\biggl(\sum_{j\ge 1}jp_j f^{j-1}(x)\biggr)G_y'(1,x),\\
G^{''}_y(1,x)&=\biggl(2\sum_{j\ge 1}jp_jG^{j-1}(y,x)G'_y(y,x)\\
&+y\sum_{j\ge 1}jp_j\bigl[(j-1)G^{j-2}(y,x)(G'_y(y,x))^2+G^{j-1}(y,x)G^{''}_y(y,x)\bigr]\biggr)\Big|_{y=1}\\
&=2G'_y(1,x)\sum_{j\ge 1}jp_j f^{j-1}(x)\\
&+(G'_y(1,x))^2\sum_{j\ge 1}j(j-1)p_jf^{j-2}(x)+G^{''}_y(1,x)\sum_{j\ge 1}jp_jf^{j-1}(x).
\end{align*}
Therefore
\begin{equation}\label{7}
\begin{aligned}
G'_y(1,x)&=\tfrac{f(x)}{1-\sum_jjp_jf^{j-1}(x)},\\
G^{''}_y(1,x)&=\tfrac{\sum_j j(j-1)p_jf^j(x)}{\biggl(1-\sum_jjp_jf^{j-1}(x)\biggr)^3}+\tfrac{2\biggl(\sum_jjp_jf^j(x)\biggr)}
{\biggl(1-\sum_jjp_jf^{j-1}(x)\biggr)^2}.
\end{aligned}
\end{equation}
So, combining \eqref{6}, \eqref{7}, and \eqref{1.1}, we conclude that
\begin{equation}\label{8}
\begin{aligned}
\Bbb E[V\Bbb I(L=\ell)]&=[x^{\ell}]\tfrac{f(x)f'(x)}{p_0},\\
\Bbb E[V(V-1)\Bbb I(L=\ell)]&=\tfrac{1}{p_0^2}\cdot [x^{\ell}]\bigl(f^{''}(x)f^2(x)\bigr),\\
&\quad+2\cdot[x^{\ell}]\bigl(\tfrac{(f'(x))^2f(x)}{p_0^2}-\tfrac{f'(x)f(x)}{p_0}\bigr).
\end{aligned}
\end{equation} 
Let us evaluate these expectations. Time for some complex analysis. Using  Weierstrass preparation theorem,
(see Ebeling \cite{Ebe},  Krantz and Parks \cite{Kra}) we proved In \cite{Pit1}  that $f(z)=\sum_{k\ge 1}z^k\, \Bbb P(L=k)$, $z\in \Bbb C$, $|z|<1$, admits an analytic extension $F(z)$ to an open disc $D$ centered at the origin, of radius $\rho>1$, {\it minus\/} a cut $[1,\rho)$ such that, for $z\in D\setminus [1,\rho)$, $z\to 1$,
\begin{equation}\label{9}
F(z)=1-\ga(1-z)^{1/2}+\sum_{s>1}\ga_s(1-z)^{s/2},\,\,\ga:=(2p_0)^{1/2}\sigma^{-1}.
\end{equation}
Here $\xi^{1/2}:=|\xi|^{1/2}\exp(i\text{Arg}(\xi)/2)$, $\text{Arg}(\xi)\in (-\pi,\pi)$. So, ever so slightly {\it above\/} the cut $\{z=x, x\in [1,\rho)\}$, we have $(1-z)^{1/2}=|1-z|^{1/2}e^{-i\pi/2}=-i|1-x|^{1/2}$. (We encourage the interested reader to check Bender \cite{Ben}, Canfield \cite{Can}, Meir and Moon \cite{MeiMoo}, \cite{MeiMoo1} and
Flajolet and Sedgewick \cite{FlaSed},  for a remarkable story of how classic complex analysis techniques found their way into analytic combinatorics.)

Since $\tfrac{F(z)F'(z)}{z^{\ell+1}}$ is integrable 
on the circular contour $\mathcal C_r=\{z= re^{i\theta},\,\theta\in (-\pi,\pi]\}$, $r=1$, we use the Cauchy integral theorem to write 
\begin{multline*}
[x^{\ell}]\bigl(f(x)f'(x)\bigr)=\tfrac{1}{2\pi i}\oint_{\mathcal C_1}\tfrac{F(z)F'(z)}{z^{\ell+1}}\,dz\\
=\tfrac{1}{2\pi i}\oint_{\mathcal C_{\rho}}\tfrac{F(z)F'(z)}{z^{\ell+1}}+\tfrac{1}{2\pi i}\int_1^{\rho}\tfrac{2i\,\text{Im}(F(x)F'(x))}{x^{\ell+1}}\,dx.
\end{multline*}
For the second line we replaced $\mathcal C_1$ with the {\it limit\/} contour.  It consists of the directed circle $z=\rho e^{i\a}$, $0\le a\le 2\pi$ with the single point $z=\rho$ pinched out, and a detour part formed by two opposite-directed line segments, one from $z=\rho e^{i(2\pi-0)}$ to $z=e^{i(2\pi-0)}$, and another from $z^{i(+0)}$ to $z=\rho^{i(+0)}$. The $2i\,\text{Im}(F(x)F'(x))$ appears because the values of $F(x)F'(x)$ just above and just below $[1,\rho)$ are complex-conjugate, so the real parts cancel each other, and $\text{Im}(F(x)F'(x))$ comes from $z$'s approaching $x\in [1,\rho)$ from above.
The first integral at the bottom is of order $O(\rho^{-\ell})$. Further, using \eqref{9} and the formula 3.191(2) from Gradschteyn and Ryzik \cite{GrR}, we get
\begin{multline*}
\tfrac{1}{2\pi i}\int_1^{\rho}\tfrac{2i\,\text{Im}(F(x)F'(x))}{x^{\ell+1}}\,dx=\tfrac{\ga}{2\pi}\int_1^{\rho}\tfrac{(x-1)^{-1/2}+O(1)}
{x^{\ell+1}}\,dx\\
=\tfrac{\ga}{2\pi}\int_1^{\infty}\tfrac{(x-1)^{-1/2}}{x^{\ell+1}}\,dx+O(1/\ell)=\tfrac{\ga (2\ell-1)!!}{2^{\ell+1}\ell!}+O(1/\ell);\\
\end{multline*}
the explicit term is of order $\ell^{-1/2}$ exactly. Therefore
\begin{equation}\label{10}
\Bbb E[V\Bbb I(L=\ell)]=\tfrac{\ga (2\ell-1)!!}{p_02^{\ell+1}\ell!}+O(1/\ell).
\end{equation}
Now, by \eqref{1.05}, we have $\Bbb P(L=\ell)=\tfrac{\ga (2\ell-3)!!}{2^{\ell} \ell!}+O(\ell^{-2})$. This and \eqref{10} imply
that
\begin{equation}\label{11}
\Bbb E[V|L=\ell]=\tfrac{\Bbb E[V\Bbb I(L=\ell)]}{\Bbb P(L=\ell)]}=\tfrac{\ell}{p_0}+O(\ell^{1/2}).
\end{equation}
Turn to the second identity in \eqref{8}. Similarly to the above computation, we obtain
\begin{align}
&\tfrac{1}{p_0^2}\cdot [x^{\ell}]\bigl(f^{''}(x)f^2(x)\bigr)=\tfrac{\ga(2\ell+1)!!}{(2p_0)^22^{\ell}\ell!}+O(1),\notag\\
&2\cdot[x^{\ell}]\bigl(\tfrac{(f'(x))^2f(x)}{p_0^2}-\tfrac{f'(x)f(x)}{p_0}\bigr)=O(1),\notag\\
\Longrightarrow\,& \Bbb E[V(V-1)\Bbb I(L=\ell)]=\tfrac{\ga(2\ell+1)!!}{(2p_0)^22^{\ell}\ell!}+O(1),\notag\\
\Longrightarrow\,& \Bbb E[V^2\Bbb I(L=\ell)]=\tfrac{\ga(2\ell+1)!!}{(2p_0)^22^{\ell}\ell!}+O(1),\notag\\
\Longrightarrow\,& \Bbb E[V^2|L=\ell]=\bigl(\tfrac{\ell}{p_0}\bigr)^2+O(\ell^{3/2}).\label{12}
\end{align}
Putting \eqref{11} and \eqref{12}, we conclude that
\begin{equation}\label{13}
\begin{aligned}
\text{var}(V|L=\ell)&:=\Bbb E\bigl[(V-\Bbb E[V|L=\ell])^2|L=\ell\bigr]\\
&=\Bbb E[V^2|L=\ell]-\Bbb E^2[V|L=\ell]=O(\ell^{3/2}).
\end{aligned}
\end{equation}
Recall that $V$ is the number of vertices in the terminal tree $T$, and $X(t)$ is the number of vertices whose descendant tree has $t$ leaves. Thus $\{\tfrac{X(t)}{V}\}_{t\ge 1}$ is the distribution of $Y$, the number of leaves in the total descendant subtree rooted at the uniformly random vertex, {\it conditioned\/} on $T$. We use \eqref{13}, {\it and\/} $X(t)\le \tfrac{\ell}{t}$, which holds since no two subtrees each with $t$ leaves intersect,  to estimate
\begin{multline*}
\Bbb E\biggl[\Bigl(\tfrac{X(t)}{V}-\tfrac{X(t)}{\Bbb E[V|L=\ell]}\Bigr)^2\,\Big| L=\ell\biggr]
\le \Bbb E\Bigl[\bigl(\tfrac{X(t)}{\ell^2}\bigr)^2\text{var}(V|L=\ell)\Bigr]=O(\ell^{-1/2}t^{-2}). \\
\end{multline*}
So, by Cauchy-Schwartz inequality, we have: uniformly for all $t\ge 1$,
\[
\Bbb E\bigl[\tfrac{X(t)}{V}|L=\ell] - \tfrac{\Bbb E[X(t)|L=\ell]}{E[V|L=\ell]}=O(\ell^{-1/4}t^{-1}).
\]
or equivalently
\begin{equation*}
\Bbb P(\mathcal L=t|L=\ell)= \tfrac{\Bbb E[X(t)|L=\ell]}{E[V|L=\ell]}+O(\ell^{-1/4}t^{-1}).
\end{equation*}
Combining this formula with Lemma \ref{lem1} and \eqref{11}, we obtain 
\begin{equation}\label{14}
\Bbb P(\mathcal L=t|L=\ell)=(1-p_1)\tfrac{(\ell-t+1)\,\Bbb P(L=\ell-t+1)\,\Bbb P(L=t)}{\ell\, \Bbb P(L=\ell)}+O(\ell^{-1/4}t^{-1}),
\end{equation}
uniformly for $t\le\ell$. Applying the formula \eqref{1.05} to \eqref{14}, we conclude
\[
\lim_{\ell\to\infty}\Bbb P(\mathcal L=t|L=\ell)=(1-p_1)\Bbb P(A_t).
\]
Since $\sum_{t\ge 1}\Bbb P(A_t)=1$, we see that the sequence of the distributions $\{\Bbb P(Y=t|L=\ell)\}_{1\le t\le \ell}$, $(\ell\ge 1)$, is not
tight, unless (as in the binary case) $\,p_1=0$. In words, $p_1$ is the limiting deficit of the leaf-set size distribution for the random subtree. 

Let us make this claim more precise.  Introduce $\tau=\tfrac{\ell^{1/2}}{\log^2\ell}$. For $t\le\tau$, it follows from \eqref{1.05} that 
\[
\tfrac{(\ell-t+1)\,\Bbb P(L=\ell-t+1)}{\ell\, \Bbb P(L=\ell)}=1+O(\tau\ell^{-1}+\ell^{-1/2})=1+O(\ell^{-1/2}),
\]
implying that 
\begin{multline*}
(1-p_1)\sum_{t\le\tau}\tfrac{(\ell-t+1)\,\Bbb P(L=\ell-t+1)\,\Bbb P(L=t)}{\ell\, \Bbb P(L=\ell)}
=(1-p_1)\sum_{t\le\tau}\Bbb P(L=t)+O(\ell^{-1/2})\\
=1-p_1 +O(\Bbb P(L>\tau))+O(\ell^{-1/2})
=1-p_1+O\bigl(\ell^{-1/4}\log\ell\bigr).
\end{multline*}
Therefore
\begin{multline*}
\Bbb P(\mathcal L>\tau|L=\ell)=1-(1-p_1)\sum_{t\le \tau}\tfrac{(\ell-t+1)\,\Bbb P(L=\ell-t+1)\,\Bbb P(L=t)}{\ell\, \Bbb P(L=\ell)}\\
+O(\ell^{-1/4}\log\ell)
=p_1+O(\ell^{-1/4}\log\ell).
\end{multline*}
This completes the proof of Theorem \ref{thm2}.

{\bf Acknowledgment.\/} I am grateful to David Aldous for helping me to put this work into a proper perspective.

\end{document}